\documentstyle[11pt]{article}

\renewcommand{\bibitem}{\item[]}
\newenvironment{biblist}
   {\begin{list}{}
    {\setlength{\leftmargin}{4ex}   
     \setlength{\itemindent}{-4ex}  
     \setlength{\itemsep}{0ex}      
     \setlength{\parsep}{0ex}       
   }}{\end{list}}

\newtheorem{prfT}{Proof of Theorem}[section]

\newtheorem{prfL}{Proof of Lemma}[section]

\newtheorem{prfC}{Proof of Corollary}[section]

\newtheorem{eg}{Example}[section]
\newenvironment{example}{\begin{eg} \rm}{\hfill $\Box$ \end{eg}}

\def \mmax {\mathop{max}\limits}

\begin{document}

\begin{center}{\Large{\textbf{Respondent privacy and estimation efficiency in randomized
response surveys for discrete-valued sensitive variables }}}
\vskip.3in
 Mausumi Bose \vskip.1in \emph{Indian Statistical
Institute, Kolkata 700108, India}
\end{center}

\vskip.2in \noindent \textbf{Abstract}: In some socio-economic
surveys, data are collected on sensitive or stigmatizing issues such
as tax evasion, criminal conviction, drug use, etc. In such surveys,
direct questioning of respondents is not of much use and the
randomized response technique is used instead.  A few researchers
have studied the issue of privacy protection or respondent jeopardy
for surveys on dichotomous populations, where the objective is to
estimate the proportion of persons bearing the sensitive trait.
However, not much is yet known about respondent protection when the
variable under study takes discrete numerical values and the
objective of the survey is to estimate the population mean of this
variable. In this article we study this issue.  We first propose a
randomization device for this situation  and give the corresponding
estimation procedure. We next propose a measure of privacy and show
that given a certain stipulated level of this privacy measure, we
can determine the parameter of the randomization device so as to
maximize the efficiency of estimation, while guaranteeing the
desired level of privacy protection.  In particular, our study also
covers the case of polychotomous populations and we can estimate the
proportions of individuals belonging to the different classes.
Consequently, results for dichotomous populations follow as
corollaries.

 \vskip.1in \noindent
\textbf{Keywords:} Jeopardy measure, numerical stigmatizing
variable, revealing probability.

\section{Introduction}\label{INT}
The randomized response technique is a useful method for collecting
data  on variables which are considered sensitive, incriminating or
stigmatizing for the respondents. Examples of such situations are
common in socio-economic surveys, for instance, we may need to
collect data on tax evasion, alcohol addiction, illegal drug use,
criminal behaviour or past criminal convictions.
 In such surveys, direct questions are not useful as the
respondents will either refuse to answer embarrassing questions or,
even if they do,  may give false answers. In a  randomized response
model,  the respondents use a randomization device
  to generate a randomized response and  the parameter under study can be estimated from these
responses. So, the respondent is not required to disclose his true
response and it is expected that this will lead to better
participation in the survey on sensitive issues.

Warner (1965) introduced the randomized response technique for
estimating the proportion of persons bearing a sensitive attribute
in a dichotomous population. In Warner's model, with population
categories $A$ and $A^c$,  a box with two types of cards labeled $A$
and $A^c$ (in proportion $p:1-p$) is used as the randomization
device. A respondent draws a card at random and responds `yes' or
`no' according as whether or not he
 belongs to the card type he draws.  Since then, several
researchers have extensively contributed to this area, e.g., Kuk
(1990), Ljungqvist (1993), Mangat
 (1994), Chua and Tsui (2000), Van den Hout and Van der Heijden (2002), Christofides
(2005) and many others. For details on the results available on this
technique we refer to the review paper by Chaudhuri and Mukerjee
(1987) and books by Chaudhuri and Mukerjee (1988) and Chaudhuri
(2011).

Lanke (1976) and Leysieffer and Warner (1976)  initiated the study
of efficiency versus privacy protection in randomized response
surveys where the population is divided into two complementary
sensitive groups, $A$ and $A^c$,  and the objective is to estimate
the proportions of persons belonging to these groups. They
 suggested measures
of jeopardy based on the `revealing probabilities', i.e., the
posterior probabilities of a respondent belonging to groups $A$ and
$A^c$ given his randomized response. Since then,
 this dichotomous case has been widely studied.  Loynes (1976)
 extended the jeopardy measure of Leysieffer and Warner (1976) to polychotomous populations.  Ljungqvist
(1993) gave a unified and utilitarian approach to measures of
privacy for the dichotomous case.
 Nayak and Adeshiyan (2009) proposed a measure of jeopardy for
surveys from dichotomous populations and developed an approach for
comparing the available randomization procedures. These results are
all based on samples drawn by simple random sampling with
replacement. 

All the references given above are for  sensitive variables which
are
 categorial or qualitative in nature. However, in randomized response surveys it is quite common to
have situations where the study variable $X$ is quantitative, e.g.
in studies on the number of criminal convictions of a person, the
number of induced abortions, the number of months spent in a
correction centre, the amount of undisclosed income,  etc. Anderson
(1977) studied the case of continuous sensitive variables and
considered the amount of information provided by the randomized
responses. For ensuring more privacy he recommended that the
expectation of the conditional variance of $X$ given the randomized
response be made as large as possible. However, not much work seems
to have been done in studying the respondent privacy aspect for
discrete-valued sensitive variables, even though surveys are often
undertaken on such variables.

To fill this gap, in this article we focus on studying the issue of
privacy protection
 when the underlying variable under study is quantitative and discrete.
  We propose the use of a randomization device and give the associated estimation method. Then, we consider
two separate cases, one where all values of $X$ are sensitive and
another where not all values of $X$ are sensitive. For each of these
cases, we propose a  measure for protecting the privacy of the
respondents. We finally show how one can choose the randomization
device parameter in each case, so as to guarantee a certain
pre-specified level of respondent protection and then maximize the
efficiency of estimating the parameter of interest under this
constraint.  Our study also covers  qualitative sensitive variables,
i.e., cases where the population is dichotomous or polychotomous,
and allows us to estimate the proportions of individuals belonging
to each category.

In Section 2 we give some preliminaries. In Sections 3 and 4 we
consider the issues of estimation and privacy protection,
respectively. In Section 5 we obtain the randomization device
parameter which allows efficient estimation while assuring the
required level of respondent protection and illustrate with some
numerical examples. In the concluding section we show how our study
covers the case of polychotomous variables.

\section{Preliminaries}\label{prelim} Consider a population
 with $N$ individuals labeled $1, \ldots, N$.  Let $X$ denote the
sensitive variable of interest. We assume that $X$ takes a finite
number of  values $x_1, \ldots, x_m$ and without loss of generality,
we may suppose these $m$ values  to be known. For $1\leq i \leq m$,
let $\pi_i$ be the unknown population proportion of individuals for
whom $X$ equals $x_i$, i.e.,
\begin{equation}\label{pii}
{\mathrm{Prob}}(X=x_i) =  \pi_i,  \ \ 1\leq i \leq m, \ \ \mathrm{where}  \ \ \pi_i \geq 0, \ \ \sum_{i=1}^m \pi_i =1, \\
\end{equation}
 The
objective of the survey is to estimate the population mean  of $X$.
 For this, we suppose as usual (cf. Warner
(1965), Nayak and Adeshiyan (2009) and others), that a sample of $n$
individuals is drawn from the population by simple random sampling
with replacement. As for the randomization device, since we are
interested in the numerical values of $X$, we propose the use of a
device as described below.

Consider a box containing cards of $(m+1)$ types, the $i$th type of
card being marked `Report $x_i$ as your response', $1\leq i\leq m$,
while the $(m+1)$th type of card is marked: `Report your true value
of $X$ as your response.' The box has a large number of cards, say
$M$, there being $Mp$ cards of type $(m+1)$ and $M\frac{1-p}{m}$
cards of each of the types $i$, $1\leq i\leq m$, $0<p<1.$ A sampled
respondent is asked to draw a card at random from the box and then
give a truthful response according to the card drawn by him, without
disclosing the label on the card to the investigator. Thus the true
value of $X$ for the respondent is not known. The $n$ responses so
received are the data from this survey.

Let $R$ denote the randomized response variable. Clearly, with this
device, the ranges of $R$ and $X$ match. The  efficiency in
estimation and respondent protection will depend on the choice of
the value of $p$, which we call the device parameter. The above
device is such that
 with probability $p$, a
respondent will report his true value,  while with probability
$\frac{1-p}{m}$, he will report any one of the possible values $x_1,
\ldots, x_m$ chosen at random, i.e.,
\begin{eqnarray}
 {\mathrm{Prob}}(R=x_i |X=x_j) &=&\frac{1-p}{m},  \ \ 1\leq i \neq j \leq m,  \\
 {\mathrm{Prob}}(R=x_j | X=x_j) &=& p +\frac{1-p}{m}, \ \ 1\leq  j \leq
 m.
\end{eqnarray}

\section{Estimation of population mean}\label{estimation}

 The population mean and variance of $X$ are given by
$$ \mu_X =\sum_{i=1}^m x_i \pi_i  \ \mbox{ and } \ \ \sigma_X^2 =
\sum_{i=1}^m (x_i -\mu_X)^2\pi_i,
$$
 respectively. Our
objective is to estimate $\mu_X$ from  the $n$ randomized responses
collected as described in Section~\ref{prelim}. Let $w_i$ be the
sample proportion of randomized responses which equal $x_i$, \
$1\leq i\leq m.$ Hence, from (1)--(3),
 \begin{equation}\label{lambda}
{\mathrm{E}}(w_i) = {\mathrm{Prob}}(R=x_i) = p\pi_i + \frac{1-p}{m}
= \lambda_i, \ \ {\mathrm{say}}. \nonumber
\end{equation} So, an unbiased estimator of $\pi_i$ will be given by
$\hat{\pi}_i = \frac{1}{p}(w_i - \frac{1-p}{m}),$  leading to an
unbiased estimator of $\mu_X$ as
$$
\hat{\mu}_X= \sum_{i=1}^m x_i\hat{\pi}_i = \frac{1}{p}\sum_{i=1}^m
x_iw_i  - \frac{1-p}{mp}\sum_{i=1}^m x_i.
$$
 Then, on simplification using
(\ref{lambda}), and writing $\bar{X}= \frac{1}{m}\sum_{i=1}^m x_i$,
the variance of $\hat{\mu}_X$ is given by
\begin{eqnarray}\label{var}
{\mathrm{Var}}(\hat{\mu}_X) &=& \frac{1}{p^2} {\mathrm{Var}}
(\sum_{i=1}^m x_iw_i ) = \frac{1}{np^2} \left\{\sum_{i=1}^m x_i^2
\lambda_i (1-\lambda_i) -  \sum\sum_{i\neq j=1}^m
x_i x_j\lambda_i \lambda_j  \right\} \nonumber \\
  & = &\frac{1}{np^2} \left\{p\sum_{i=1}^m  x_i^2\pi_i +
  \frac{1-p}{m}\sum_{i=1}^m  x_i^2 - \left(p\mu_X + (1-p)\bar{X}\right)^2
  \right\} \nonumber \\
  & =& \frac{1}{np^2}\left\{ p \sigma_X^2 + (1-p) \frac{1}{m}\sum_{i=1}^m
(x_i-\bar{X})^2 + p(p-1) (\mu_X - \bar{X})^2
  \right\}.
\end{eqnarray}

 Our aim is to estimate $\mu_X$ keeping
{Var}$(\hat{\mu}_X)$ as small as possible. It is clear from the
expression on the right side of (\ref{var}) that
{Var}$(\hat{\mu}_X)$ is decreasing in $p$, irrespective of the
values of $\pi_1, \ldots, \pi_m$. So, this variance may be
decreased, or equivalently, the efficiency of estimation may be
increased by increasing $p$, whatever may be the proportions of the
$x_i$ values in the population.


\section{Privacy protection}\label{privacy}

To study the respondent privacy aspect for dichotomous populations,
Leysieffer and Warner
 (1976) studied the case where both $A$ and $A^c$ are sensitive categories while Lanke (1975) also considered the case where
 only $A$ is sensitive and
there is no jeopardy in a `no' answer to the sensitive question. For
 polychotomous populations, Loynes (1976) studied two cases, one where
all categories are stigmatizing and another where one of the
categories is not stigmatizing. In line with these,
 we too consider the
privacy issue for two situations, one where all the $m$ values of
$X$ are stigmatizing and another where not all  values of $X$ are
 stigmatizing.  Both these situations  commonly arise in
practice and we  require separate privacy protection measures for
them.

For a randomly chosen respondent from the population, the `true'
probability that the value of  $X$ for this respondent equals $x_i$
is given by Prob$(X=x_i)$. On the other hand, when this respondent
gives a randomized response, say $x_j$, then the probability that
the value of $X$ for this respondent equals $x_i$ is now given by
the conditional probability Prob$(X=x_i|R=x_j)$, or the  `revealing'
probability.

\subsection{All values of $X$ are stigmatizing}\label{all}
Suppose all the values $x_1, \ldots, x_m$ are stigmatizing. In this
case, a respondent would feel comfortable in participating in the
survey if the perception of his having a value $X=x_i$ is not much
altered after knowing his randomized response, for all $1\leq i\leq
m$. This would require that his true and revealing probabilities be
sufficiently close. Starting from this basic premise we define
\begin{equation}\label{alphaij}
\alpha_{ij} = |{\mathrm Prob}(X=x_i|R=x_j) -{\mathrm Prob}(X=x_i)|
\end{equation}
and since each respondent would want $\alpha_{ij}$ to be as small as
possible for all $1\leq i,j \leq m$, as a measure of privacy
protection we propose the following measure:
\begin{equation}\label{alpha}
 \alpha = \mmax_{1\leq i,j \leq m}
\alpha_{ij}.\end{equation}

 A randomization device with a privacy protection
value $\alpha=\alpha_0$ would guarantee that the discrepancies
between the true  and revealing probabilities will be at most
$\alpha_0$ for all respondents, irrespective of their true values.
Thus a device which results in a lower value of $\alpha$ gives a
higher level of privacy protection  than one with a higher value of
$\alpha$.

Suppose the scientist planning a certain survey  would like to keep
the privacy protection available to respondents above a certain
threshold, i.e., would like to achieve $\alpha \leq \xi$, where
$\xi$ is a pre-assigned  quantity, $0< \xi < 1$. Moreover, this
bound on $\alpha$ should hold irrespective of the unknown values of
$\pi_1, \ldots, \pi_m.$ The following theorem shows how the
 device parameter can be chosen to achieve this.

\noindent {\bf Theorem 1.} {\it For $\alpha$ as in (\ref{alpha}) and
a preassigned $\xi$, where  $0< \xi < 1$,  \ $\alpha \leq \xi $ will
hold, irrespective of the values of $\pi_1, \ldots, \pi_m, $ if and
only if $p \leq p_0$,  where }
\begin{equation}\label{p01}
 p_0 = \frac{1}{1+ \frac{m}{\xi}(\frac{1-\xi}{2})^2}.
\end{equation}

\noindent {\bf Proof.} From (1)-(3), using Bayes' Theorem it follows
that for $1\leq i,j,\leq m$,
 \begin{equation}\label{revealing}
 {\mathrm Prob} (X=x_i|R=x_j) = \frac{(p\delta_{ij} +
\frac{1-p}{m})\pi_i}{\sum_{u=1}^m(p\delta_{ju} +
\frac{1-p}{m})\pi_u}= \frac{(p\delta_{ij} +
\frac{1-p}{m})\pi_i}{p\pi_j + \frac{1-p}{m}},
\end{equation}
where $\delta_{ij}$ is Kronecker Delta. Hence  from (\ref{alphaij})
it follows that
$
 \alpha_{ij} = \frac{p\pi_i|\pi_j-\delta_{ij}|}{p\pi_j + \frac{1-p}{m}}$
and for any $i\neq j$,
\begin{equation}
 \alpha_{ij} =\frac{p\pi_i\pi_j}{p\pi_j + \frac{1-p}{m}}\leq \frac{p(1-\pi_j)\pi_j}{p\pi_j +
 \frac{1-p}{m}} = \alpha_{jj}, \nonumber
\end{equation}
 as $\pi_i + \pi_j
\leq 1$ for all $i,j$. Thus $ \alpha = \mmax_{1\leq j\leq m}
\alpha_{jj} = \mmax_{1\leq j\leq m} \frac{\pi_j(1-\pi_j)}{\pi_j +
 \frac{1-p}{mp}}.
$ Hence, $\alpha \leq \xi$
 if  and only if
\begin{equation}\label{inequality}
\pi_j(1-\pi_j) -\xi\pi_j \leq \frac{\xi(1-p)}{mp} \ \  {\mathrm for
\ \ all} \ \ 1\leq j \leq m.
\end{equation}
First suppose $p\leq p_0$. Then for $1\leq j \leq m$,
\begin{eqnarray}
 \pi_j(1-\pi_j) -\xi\pi_j  &=&  \left(\frac{1-\xi}{2}\right)^2 - \left(\frac{1-\xi}{2} - \pi_j \right)^2 \nonumber \\
  &\leq & \left(\frac{1-\xi}{2}\right)^2 = \frac{\xi(1-p_0)}{mp_0}, \ \ {\mathrm using } \ \ (\ref{p01})  \nonumber \\
   &\leq & \frac{\xi(1-p)}{mp}. \nonumber
\end{eqnarray}
Thus the inequalities in (\ref{inequality}) hold, or equivalently
$\alpha \leq \xi$, irrespective of the values of $\pi_1, \ldots,
\pi_m.$

To prove the converse,  suppose  $\alpha \leq \xi$, or equivalently,
the inequalities in (\ref{inequality}) hold, irrespective of the
values of $\pi_1, \ldots, \pi_m.$ Then, for $\pi_1=\frac{1-\xi}{2},
\pi_2 = \frac{1+\xi}{2}, \pi_3 = \ldots = \pi_m=0$,  in particular,
these inequalities will also hold. So, for this choice of $\pi_j$
values in (\ref{inequality}) with $j=1$,  we have
\begin{eqnarray}
  \left(\frac{1-\xi}{2}\right)\left(\frac{1+\xi} {2}\right) - \xi\left(\frac{1-\xi}{2}\right)& \leq &
  \frac{\xi(1-p)}{mp} \nonumber \\
{ \mathrm i.e., } \ \ \left(\frac{1-\xi}{2} \right)^2   & \leq &
\frac{\xi(1-p)}{mp}. \nonumber
\end{eqnarray}
So from (\ref{p01}),
$p\leq p_0$. Hence theorem.
\hfill$\Box$

\noindent \textbf{Remark 1.} It is clear from (\ref{p01}) that in
order to maintain the same level of protection, the value of $p_0$
monotonically decreases with the number of possible values of $X$.
Again, for a given number of possible values of $X$, $p_0$
monotonically increases with $\xi$. We may reiterate that these
values of $p$ do not depend on how the values of $X$ are distributed
in the population.

\subsection{Not all  values of $X$ are  stigmatizing}

In many surveys it may so happen that not all values of $X$ are
sensitive or stigmatizing. For instance, in a survey for estimating
the average number of criminal convictions of persons in a certain
population, the value $X=0$ is not stigmatizing but any value of $X
\geq 1$ could well be stigmatizing. Similarly, for a survey for
estimating the average of the  number (X)  of induced abortions, the
values  $X=0$ or $X=1$ might not be considered  as stigmatizing
values while other larger values might be considered stigmatizing by
the respondents.

To study the respondents' privacy protection for such surveys, we
present here the simpler case where only one of the values of $X$,
say $x_1$, is not stigmatizing, while values $x_2, \ldots, x_m$ are
considered stigmatizing. We develop the
 protection measure for this case in detail.
Later we remark that the results obtained for this case may be
easily extended to the case where $X$ has more than one
non-stigmatizing values.

As before, the data collection and estimation proceeds as in
Sections~\ref{prelim} and \ref{estimation}. To study the respondent
protection  we note that since the value $x_1$ is non-stigmatizing,
 respondents  will feel comfortable with a randomization device
for which
 the `revealing' probability of their having  a true value
$x_1$  will be large.  So,  we propose the following measure of
privacy:
\begin{equation}\label{beta}
 \beta = \min_{1\leq j\leq m} P(X=x_1|R=x_j)
 = \min_{1\leq j\leq m} \frac{(p\delta_{1j} + \frac{1-p}{m})\pi_1}{p\pi_j +
\frac{1-p}{m}},
\end{equation}
 on simplification using (\ref{revealing}).
A device with a privacy protection value $\beta$ will guarantee that
all respondents are perceived to have $X=x_1$ with probability at
least $\beta$. So, a device leading to a larger value of $\beta$
will ensure greater privacy to respondents than one with a smaller
$\beta$.

 Let $\xi$, $0<\xi<1$, denote a preassigned level of
 respondents' privacy. Then in order to achieve this level
of protection we require that $\beta \geq \xi$, irrespective of the
values of $\pi_1, \ldots, \pi_m$. Thus we should have
$$(p\delta_{1j}+\frac{1-p}{m})\pi_1 \geq \xi (p\pi_j +
\frac{1-p}{m}), \ \ \ 1\leq j\leq m,$$  or equivalently,    the
following inequalities should hold:
\begin{eqnarray}\label{ineq1}
[p(1-\xi) + \frac{1-p}{m}]\pi_1 &\geq &\frac{\xi(1-p)}{m} \\ {\rm
 and}  \ \ \ \frac{1-p}{m}\pi_1 -\xi p\pi_j & \geq &\frac{\xi(1-p)}{m}, \ \
2\leq j \leq m.
\end{eqnarray}
Clearly, no $p$ can satisfy (13) irrespective of $\pi_1, \ldots,
\pi_m$ for any given $\xi$ since (13) fails as $\pi_1 \rightarrow
0.$ So we assume that $\pi_1 > 0$ and we also assume some prior
knowledge about a lower bound on $\pi_1$. This assumption is quite
realistic because  in most populations there will be an appreciable
number of  persons with a non-stigmatizing variable value and hence,
a lower bound to the proportion of such stigma-free persons in the
population will be available.

Thus, suppose we have prior knowledge that $\pi_1 \geq c$. We work
with $\xi<c$. This  is again realistic because if the only knowledge
about $\pi_1$ is that $\pi_1\geq c$, it is impractical to demand
that $P(X=x_1|R=x_j) \geq \xi (\geq c) $ for all $j$.  Now,  the
following theorem gives the value of the device parameter $p$ which
will guarantee the desired level of respondent protection $\xi$.

\noindent \textbf{Theorem 2.} {\it Let $\beta$ be as in (\ref{beta})
and $\pi_1 \geq c$ for some known $c$. Then given a preassigned
$\xi$, where $0< \xi < c $,
 \ $\beta \geq \xi $ will hold,
irrespective of the values of $\pi_1, \ldots, \pi_m, $ if and only
if $p \leq p_0$, where }
\begin{equation}\label{p02}
 p_0 = \frac{\frac{c-\xi}{m}}{\frac{c-\xi}{m} + \xi(1-c)}.
\end{equation}

\noindent {\bf Proof.} Since $\pi_1 \geq c$, it is clear that $\pi_j
\leq 1-c$ for $2\leq j\leq m$ and  we have
\begin{eqnarray}
[p(1-\xi) + \frac{1-p}{m}]\pi_1 &\geq & [p(1-\xi) + \frac{1-p}{m}]c
\nonumber\\ {\rm and }  \ \ \ \frac{1-p}{m}\pi_1 - \xi p\pi_j & \geq
& \frac{1-p}{m} c - \xi p(1-c), \ \ 2\leq j \leq m. \nonumber
\end{eqnarray}
As a result, (\ref{ineq1}) and (14) will hold, irrespective of the
true values of $\pi_1(\geq c), \pi_2, \ldots, \pi_m$ iff
\begin{eqnarray}\label{ineq3}
[p(1-\xi) + \frac{1-p}{m}]c &\geq & \xi \frac{1-p}{m}\\ {\rm and } \
\ \  \frac{1-p}{m}c -\xi p(1-c) & \geq & \xi \frac{1-p}{m}
\end{eqnarray}
hold. Now, (\ref{ineq3}) reduces to
$$(p+\frac{1-p}{m})c \geq \xi (cp + \frac{1-p}{m})$$ which will
always hold for every $p$ since $\xi (cp+\frac{1-p}{m}) \leq
\xi(p+\frac{1-p}{m}) < c(p+\frac{1-p}{m})$ as $\xi<c$ and
$p+\frac{1-p}{m} >0.$  So, it is enough to only consider (17). Note
that \begin{eqnarray} (17)  \Leftrightarrow  \frac{c-c p}{m} - \xi
p(1-c) &
\geq & \frac{\xi-\xi p}{m} \nonumber \\
\Leftrightarrow p & \leq & = \frac{\frac{c-\xi}{m}}{\frac{c-\xi}{m}
+ \xi(1-c)} = p_0, \nonumber
\end{eqnarray}thus proving the theorem. \hfill$\Box$

\noindent \textbf{Remark.} The above discussion can be extended to
include the more general case where $X$ has $t$ non-stigmatizing
values $x_1, \ldots, x_t$, say, while its remaining $m-t$ values are
stigmatizing, $1<t<m.$ In that case too, it can be shown that $p_0$
takes the form as in Theorem 2, but now with $$\beta=\\min_{1\leq j
\leq m} P(X=x_1 {\mbox { or }} x_2 {\mbox { or }}\ldots x_t|R=x_j)
{\mbox {  and  }} \pi_1+\ldots + \pi_t \geq c {\mbox {  with
}}\xi<c.$$

\section{Privacy protection together with efficiency in
estimation}\label{privacyandefficiency} We now consider the issue of
efficiency in estimation together with
 privacy protection in randomized response surveys. It was seen from (\ref{var}) that,  irrespective of the
values of $\pi_1, \ldots, \pi_m$, the efficiency of estimation may
be increased by increasing $p$. On the other hand, for a given $\xi$
and irrespective of the values of $\pi_1, \ldots, \pi_m$, Theorems 1
and 2 show that  a protection of $\alpha \leq \xi$ or $\beta \geq
\xi$ may be guaranteed iff $p\leq p_0$, where $p_0$ is as in
(\ref{p01}) or (\ref{p02}), respectively.   So, the best choice of
$p$ with regard to maximizing the efficiency of estimation of
$\mu_X$, subject to the stipulated level of privacy protection
$\xi$, is $p=p_0$.   The following examples illustrate this.
\begin{example}  Let $X$ take four  values which are all sensitive.  Suppose  $\xi =0.1$ Then by Theorem 1,  $p_0 =
 0.1099.$ So, if we use a randomization
device with  $p=0.1099$ then the efficiency of estimation can be
maximized while  guaranteeing that the maximum discrepancy between
the true probability and the revealing probability of all
respondents will be at most 0.1.
\end{example}
 The following table gives the $p_0$ values in (\ref{p01}) for some
choices of $\xi$ and $m$.
\begin{center}
\begin{tabular}{|c|c|c|c|c|c|c|c|c|c|c|}
\hline
    $m$ & $\xi$ & $p_0$      &  & $m$ & $\xi$ & $p_0$      &  & $m$ & $\xi$ & $p_0$ \\ \hline
    3 & 0.1 & 0.1413 &  & 4 & 0.1 & 0.1099 &  & 5 & 0.1 & 0.0899  \\
    3 & 0.2 & 0.2941 &  & 4 & 0.2 & 0.2381 &  & 5 & 0.2 &  0.2000\\
    3 & 0.3 & 0.4494 &  & 4 & 0.3 & 0.3797 &  & 5 & 0.3 & 0.3288 \\
    3 & 0.4 & 0.5970 &  & 4 & 0.4 & 0.5263 &  & 5 & 0.4 & 0.4706 \\
    \hline
  \end{tabular}
\end{center}
\begin{example}  Let $X$ take one nonsensitive value and two sensitive
values. Suppose it can be assumed that at least 15\% of the
 individuals in the population possess the nonsensitive value and suppose it
is stipulated that $\xi= 0.10$. Then by Theorem 2, $p_0=
0.1639$. So, if we use a device with $p=0.1639$ then estimation
efficiency will be maximum while guaranteeing that all respondents
will have at least a 10\% probability of being revealed as belonging
to the non-stigmatizing class.
\end{example}
\section{Estimation of population proportions} As mentioned in
Section~\ref{INT}, several researchers have estimated the
proportions of individuals belonging to the two categories in
dichotomous populations, while Loynes (1976) extended this to
estimating the different proportions in a
 polychotomous population.
In our case where $X$ takes $m$ numerical values, we may also
readily estimate the population proportions $\pi_1, \ldots, \pi_m$
from the responses collected as in Section~\ref{prelim} and again
use the measures of privacy as given in (\ref{alpha}) and
(\ref{beta}) to achieve the stipulated level of privacy protection.

As seen in Section~\ref {estimation}, an unbiased estimate of
$\pi_i$ is
$$\hat{\pi}_i = \frac{1}{p}(w_i -\frac{1-p}{m}), \ \ \ \ 1\leq i\leq
m.$$ Suppose, in the spirit of $A-$optimality commonly used in
optimal design theory,  we would like to minimize the average
variance of these estimates. For this, we can show that
\begin{eqnarray}\label{varpi}
\sum_{i=1}^m Var (\hat{\pi_i})  = \frac{1}{np^2} \sum_{i=1}^m
\lambda_i (1-\lambda_i)
 = \frac{1}{n}\left\{ \frac{1}{p^2} -  \sum_{i=1}^m \pi_i^2 + \frac{1}{m}(\frac{1}{p^2} -1)\right\},
\end{eqnarray}
on simplification, using (\ref{lambda}). Clearly, (\ref{varpi}) is
decreasing in $p$, irrespective of the true values of $\pi_1,
\ldots, \pi_m.$ So as in the case of estimating the mean, here too,
given some $\xi$,  subject to the constraint on protection of
privacy, the best choice for $p$ for minimizing the average variance
of the estimates of the proportions, is $p=p_0$, with $p_0$ being
given by ({\ref{p01}) or ({\ref{p02}), as the case may be. The
popular case of dichotomous populations follow by taking $m=2$ in
the above.

\end{document}